\documentclass{amsart}
\usepackage{amssymb}

\newtheorem{theorem}{Theorem}[section]
\newtheorem{corollary}[theorem]{Corollary}

\newtheorem{lemma}[theorem]{Lemma}
\newtheorem{proposition}[theorem]{Proposition}

\numberwithin{equation}{section}

\def\limfunc#1{\mathop{\rm #1}}%

\begin{document}
\title[Test groups for Whitehead groups]{Test groups for Whitehead groups}
\author{Paul C. Eklof}
\address[Eklof]{Mathematics Dept, UCI\\
Irvine, CA 92697-3875} \email{peklof@math.uci.edu}
\thanks{}
\author{L\'{a}szl\'{o}  Fuchs}
\address[Fuchs]{Department of Mathematics, Tulane University, New Orleans, LA 70118\\}
 \email{fuchs@tulane.edu}
\thanks{}
\author{Saharon Shelah}
\address[Shelah]{Institute of Mathematics, Hebrew University\\
Jerusalem 91904, Israel}
\email{shelah@math.huji.ac.il}
\thanks{ The third author would like to thank the United States-Israel Binational
      Science Foundation for their support. Publication 879.
      }
\date{\today}

\begin{abstract}
We consider the question of when the dual of a Whitehead group is a
test group for Whitehead groups. This turns out to be equivalent to
the question of when the tensor product of two Whitehead groups is
Whitehead. We investigate what happens in different models of set
theory.
\end{abstract}

\maketitle

\section{Introduction}

All groups in this note are abelian. A {\it Whitehead group}, or {\it W-group} for short, is defined to be an abelian group $A$ such that $\limfunc{%
Ext}(A,\mathbb{Z})=0$. We are looking for groups $C$ other than
$\mathbb{Z}$ such
that a group $A$ is a W-group if and only if $\limfunc{Ext}(A,C)=0$; such a $%
C$ will be called \textit{a test group for Whitehead groups}, or a \textit{%
W-test group} for short\textit{.} Notice that if $C$ is a non-zero
separable torsion-free group, then $\limfunc{Ext}(A,C)=0$ implies
$A$ is a W-group (since $\mathbb{Z}$ is a summand of $C$), but the
converse may not hold, that is, $\limfunc{Ext}(A,C)$ may be non-zero
for some W-group $A$.

Among the separable torsion-free groups are the dual groups, where
by a dual group we
mean one of the form $\limfunc{Hom}(B,\mathbb{Z})$ for some group $B$. We call $%
\limfunc{Hom}(B,\mathbb{Z})$ the ($\mathbb{Z}$-){\it dual} of $B$ and denote it by $%
B^{*}$. We shall call a group a $W^{*}$\textit{-group} if it is the
dual of
a W-group. The principal question we will consider is whether every $W^{*}$%
-group is a W-test group. This turns out to be equivalent to a
question about tensor products of W-groups. (See the end of section
2.)

As is almost always the case with problems related to Whitehead
groups, the answer depends on the chosen model of set theory. We
have an easy affirmative answer if every W-group is free (for
example in a model of V = L); therefore we will focus on models
where there are non-free W-groups. We will exhibit models with
differing results about whether $W^{*}$-groups are W-test groups,
including information about some of the ``classical'' models where
there are non-free W-groups. In particular, we will show that the
answer to the question is independent of ZFC + GCH.

\bigskip

\section{Theorems of ZFC}

We begin with some theorems of ZFC.

\begin{theorem}
  A group $A$ is a $W^{*}$-group
if and only if it is the kernel of an epimorphism
$\mathbb{Z}^{\kappa }\rightarrow \mathbb{Z}^{\lambda }$ for some
cardinals $\kappa$ and $\lambda$.
\end{theorem}

\noindent \textsc{Proof}. If $A=B^{*}$ where $B$ is a W-group,
choose a free resolution
\[
0\rightarrow K\rightarrow F\rightarrow B\rightarrow 0
\]
where $F$ and $K$ are free groups. Taking the dual yields
\[
0\rightarrow B^{*}=A\rightarrow F^{*}\rightarrow K^{*}\rightarrow \limfunc{%
Ext}(B,\mathbb{Z})=0\text{. }
\]
Since $F^{*}$ and $K^{*}$ are products, this proves the claim in
one direction.

Conversely, if $H$ is the kernel of an epimorphism $\phi
:\mathbb{Z}^{\kappa }\rightarrow \mathbb{Z}^{\lambda }$, then we
have an exact sequence
\[
0\rightarrow H\rightarrow F^{*}\rightarrow K^{*}\rightarrow 0
\]
for certain free groups $F$ and $K$. Dualizing, we obtain an exact
sequence
\[
0\rightarrow K\rightarrow F\rightarrow B\rightarrow 0
\]
for a subgroup $B$ of $K^{*}$. When we take the dual of the last
sequence we
obtain the original epimorphism $\phi $, and conclude that its kernel is $%
B^{*}$ and $\limfunc{Ext}(B,\mathbb{Z})=0$. \qed

\begin{theorem}
For W-groups $A$ and $B$, there are natural isomorphisms
\[
\limfunc{Ext}(A,B^{*})\cong \limfunc{Ext}(B,A^{*})
\]
and
\[
\limfunc{Ext}(A,B^{*})\cong \limfunc{Ext}(A\otimes
B,\mathbb{Z})\text{.}
\]
\end{theorem}

\noindent \textsc{Proof}. Recall that, for any pair of groups $A$
and $B$, there are natural isomorphisms
\[
\limfunc{Hom}(A,B^{*})\cong \limfunc{Hom}(A\otimes B,\mathbb{Z})\cong \limfunc{%
Hom}(B\otimes A,\mathbb{Z})\cong \limfunc{Hom}(B,A^{*})\text{.}
\]

Now let $A$ and $B$ be W-groups. Using these isomorphisms as
vertical maps
and a free resolution $0\rightarrow K\rightarrow F\rightarrow A\rightarrow 0$%
, we can form a commutative diagram with exact rows as follows:
\[
\begin{array}{cccccc}
0\rightarrow  & \limfunc{Hom}(A,B^{*})\rightarrow  & \limfunc{Hom}%
(F,B^{*})\rightarrow  & \limfunc{Hom}(K,B^{*})\rightarrow  & \limfunc{Ext}%
(A,B^{*}) & \rightarrow 0 \\
& \downarrow  & \downarrow  & \downarrow  &  &  \\
0\rightarrow  & \limfunc{Hom}(A\otimes B,\mathbb{Z})\rightarrow  & \limfunc{Hom}%
(F\otimes B,\mathbb{Z})\rightarrow  & \limfunc{Hom}(K\otimes B,\mathbb{Z}%
)\rightarrow  & \limfunc{Ext}(A\otimes B,\mathbb{Z}) & \rightarrow 0 \\
& \downarrow  & \downarrow  & \downarrow  & \downarrow  &  \\
0\rightarrow  & \limfunc{Hom}(B\otimes A,\mathbb{Z})\rightarrow  & \limfunc{Hom}%
(B\otimes F,\mathbb{Z})\rightarrow  & \limfunc{Hom}(B\otimes K,\mathbb{Z}%
)\rightarrow  & \limfunc{Ext}(B\otimes A,\mathbb{Z}) & \rightarrow 0 \\
& \downarrow  & \downarrow  & \downarrow  &  &  \\
0\rightarrow  & \limfunc{Hom}(B,A^{*})\rightarrow  & \limfunc{Hom}%
(B,F^{*})\rightarrow  & \limfunc{Hom}(B,K^{*})\rightarrow  & \limfunc{Ext}%
(B,A^{*}) & \rightarrow 0
\end{array}
\]

The $0$'s at the end of the four rows are justified by the freeness of $F$%
and by the fact that $F\otimes B$, $B\otimes F$, and $B$ are
W-groups, respectively. All the horizontal and vertical maps are
natural isomorphisms. Consequently, the diagram can be
completed---preserving commutativity---by natural maps between the
two top Exts and the two bottom Exts. These maps provide the desired
natural isomorphisms. \qed

As an immediate corollary we obtain:

\begin{corollary}
\label{cor}
 Let $A$ and $B$ be W-groups. The following are equivalent:

(i) $A\otimes B$ is a W-group;

(ii) $\limfunc{Ext}(A,B^{*})=0$;

(iii) $\limfunc{Ext}(B,A^{*})=0$. \qed
\end{corollary}

The following Proposition demonstrates that the hypothesis that $A$
and $B$ are W-groups is necessary for the equivalent conditions of
Corollary \ref{cor} to hold.

\begin{proposition}
For any abelian groups $A$ and $B$, if either

(i) $A\otimes B$ is a non-zero W-group, or

(ii) $\limfunc{Ext}(A,B^{*})=0$ and $B^{*}$ is non-zero,

\noindent then $A$ is a W-group.
\end{proposition}

\noindent \textsc{Proof}. For (ii), the conclusion follows since
$B^{*}$ is separable. For (i), we use the fact (\cite[p.116]{CE})
that there is a (non-natural) isomorphism
\[
\limfunc{Ext}(A,\limfunc{Hom}(B,\mathbb{Z}))\oplus \limfunc{Hom}(A,\limfunc{Ext}%
(B,\mathbb{Z}))\cong \limfunc{Ext}(A\otimes B,\mathbb{Z})\oplus \limfunc{Hom}(%
\limfunc{Tor}(A,B),\mathbb{Z})
\]
In our case the right-hand side reduces to $\limfunc{Ext}(A\otimes
B,\mathbb{Z})$, so the hypothesis implies that
$\limfunc{Ext}(A,B^{*})=0$. It suffices then to show that $B^{*}
\neq 0$, but this follows from the fact that $A\otimes B$ is a
W-group, hence is separable, and is also assumed to be non-zero.
\qed

\bigskip

\section {An Independence Result}

In this section we will exhibit

(1) a model of ZFC + GCH such that there is a W-group $B$ of
cardinality $\aleph _{1}$ such that $B^{*}$ is not a test group for
W--groups of cardinality $\aleph _{1}$, and

(2) another model of ZFC + GCH in which there are non-free W-groups
and every W-group $B$ (of arbitrary cardinality) is a test group for
W-groups (of arbitrary cardinality).

\smallskip\

\noindent \textbf{Model (1).} For the first model, we use the model
in \cite[Theorem 0.5]{ES} in which there exists a non-reflexive
W-group. In fact, the following theorem is proved:

\begin{theorem}
It is consistent with ZFC + GCH that there is a non-free W-group $B$
of cardinality $\aleph _{1}$ such that $B^{*}$ is free.
\end{theorem}

Let $B$ be as in the theorem. Note that $B^{*}$ must be of infinite
rank, since $B$ is not free, but there is a monomorphism $: B
\rightarrow B^{**}$. Since $B$ is a non-reflexive W-group, a result
of Huber (see \cite{Hu} or
\cite[XI.2.7]{EM}) implies that $B$ is not $\aleph _{1}$-coseparable, i.e., $%
\limfunc{Ext}(B,\mathbb{Z}^{(\omega )})\neq 0$.  Thus $\limfunc{Ext}%
(B,B^{*})\neq 0$.

\medskip

\noindent \textbf{Model (2).} We shall work in a model of Ax($S$) + $%
\diamondsuit ^{*}(\omega _{1}\setminus S)$ plus

\begin{quote}
\smallskip \ ($\ddagger $) $\diamondsuit _{\kappa }(E)$ holds for every
regular cardinal $\kappa >\aleph _{1}$ and every stationary subset $E$ of $%
\kappa $.
\end{quote}

For the definition and implications of Ax($S$) + $\diamondsuit
^{*}(\omega _{1}\setminus S)$, see \cite[pp. 178-179 and Thm.
XII.2.1]{EM}. This is the ``classical'' model of ZFC + GCH in which
there are non-free W-groups of cardinality $\aleph _{1}$; in fact,
whether an $\aleph _{1}$-free group of cardinality $\aleph _{1}$ is
a W-group is determined by its Gamma invariant.

\begin{theorem}
\label{axs} Assuming $\limfunc{Ax}(S)$ $+ \diamondsuit ^{*}(\omega
_{1}\setminus S)$, if $A$ and $B$ are W-groups of cardinality $\leq
\aleph _{1}$, then $A\otimes B$ is a W-group.
\end{theorem}

\noindent \textsc{Proof}. The case when either $A$ or $B$ is
countable
(hence free) is trivial so we can assume that $A$ and $B$ have cardinality $%
\aleph _{1}$. By Theorem XII.2.1 of \cite{EM}, $\Gamma (A)\leq
\tilde{S}$ and $\Gamma (B)\leq \tilde{S}$ since $A$ and $B$ are
W-groups. It suffices
to prove that $\Gamma (A\otimes B)\leq \tilde{S}$. Fix $\omega _{1}$%
-filtrations $\{A_{\nu }:\nu <\omega _{1}\}$ and $\{B_{\nu }:\nu
<\omega _{1}\}$ of $A$ and $B$, respectively, such that
\[
\{\nu <\omega _{1}:\exists \mu >\nu \text{ s.t. }A_{\mu }/A_{\nu
}\text{ is not free}\}\subseteq S
\]
and similarly for $B$. Then $\{A_{\nu }\otimes B_{\nu }:\nu <\omega
_{1}\}$ is an $\omega _{1}$-filtration of $A\otimes B$. For each
$\mu >\nu $ there is an exact sequence
\[
0\rightarrow A_{\nu }\otimes B_{\nu }\rightarrow A_{\mu }\otimes
B_{\mu }\rightarrow (A_{\mu }\otimes (B_{\mu }/B_{\nu }))\oplus
((A_{\mu }/A_{\nu })\otimes B_{\mu })
\]
(cf. \cite[Prop. 4.3a(c), p. 25]{CE}). If $\nu \notin S$, then the
two summands on the right are free, and hence the quotient $A_{\mu
}\otimes B_{\mu }/A_{\nu }\otimes B_{\nu }$ is free. Thus we have
proved that $\Gamma (A\otimes B)\leq \tilde{S}$.\qed

\smallskip\

The following is proved by the methods of proof of Theorem 3.1 in
\cite{BFS} (see also Corollary XII.1.13 of \cite{EM}).

\smallskip \

\begin{lemma}
\label{dagger}In any model of ($\ddagger $), every W-group $A$ is
the union of
a continuous chain of subgroups $\{A_{\nu }:\nu <\sigma \}$ such that $%
A_{0}=0$ and for all $\nu <\sigma $, $A_{\nu }$ is a W-group and
$A_{\nu +1}/A_{\nu }$ has cardinality $\leq \aleph _{1}$ and is a
W-group. \qed
\end{lemma}

\medskip

The hypothesis Ax($S$) + $\diamondsuit ^{*}(\omega _{1}\setminus S)$
 implies that the second assumption of the following theorem holds. Thus, we will show that
 our
 model has the desired property if we prove the
following theorem.

\bigskip

\begin{theorem}
Assume that ($\ddagger $) holds and that the tensor product of two
W-groups of cardinality $\leq \aleph _{1}$ is again a W-group. Then
the tensor product of any two W-groups (of arbitrary cardinality) is
again a W-group.
\end{theorem}

\smallskip

\noindent \noindent \textsc{Proof}. Let $A$ and $B$ be W-groups. It
suffices to prove that $A \otimes B$ is a W-group when at least one
of the groups, say $B$, has cardinality $>\aleph _{1} $. Let
$\{B_{\nu }:\nu <\sigma \}$ be a continuous chain as in the
conclusion of Lemma \ref{dagger}. We shall first prove that if
$|A|=\aleph _{1}$, then $A\otimes B$ is a W-group. Now $\{A\otimes
B_{\nu }:\nu <\sigma \}$ is a continuous filtration of $A\otimes B$.
To show that $A\otimes B$ is a W-group, it suffices to show that for
all $\nu <\sigma $, $A\otimes B_{\nu +1}/A\otimes B_{\nu }$ is a
W-group (cf. \cite[XII.1.5]{EM}). Now, as above, there is an exact
sequence
\[
0\rightarrow A\otimes B_{\nu }\rightarrow A\otimes B_{\nu
+1}\rightarrow A\otimes (B_{\nu +1}/B_{\nu })
\]
and by hypothesis, the right-hand term is a W-group (of cardinality $%
\aleph _{1}$). But then the quotient $A\otimes B_{\nu +1}/A\otimes
B_{\nu }$ is a subgroup of a W-group and hence a W-group.

Next suppose that $|A|\grave{>}\aleph _{1}$.  Again we use the
continuous filtration $\{A\otimes B_{\nu }:\nu <\sigma \}$and the
displayed exact sequence above. Now the right-hand term $A\otimes
(B_{\nu +1}/B_{\nu })$ is a W-group by the first case (and the
symmetry of the tensor product) because it is a tensor product of
two W-groups one of which has cardinality $\aleph _{1}$. \qed

\medskip

{\bf Remark.} As mentioned above, Martin Huber  has shown that every
$\aleph_1$-coseparable group is reflexive. In Model (2) above and
also in the model of the next section, all W-groups are
$\aleph_1$-coseparable.  This raises the question of whether
(provably in ZFC) every $\aleph_1$-coseparable group, or every
reflexive group, $B$ satisfies Ext$(B, B^*)=0$, or even satisfies:
$B^*$ is a W-test group.

\bigskip

\section{Martin's Axiom}

A model of Martin's Axiom (MA) was the first model in which it was
proved (in \cite{Sh74}) that there are non-free W-groups. So it is
of interest to see what happens in that model. In fact, the
conclusion is the same as that of Theorem \ref{axs}. But we prove it
by proving the following theorem:

\begin{theorem}
(MA + $\lnot $CH) If $A$ and $B$ are W-groups of cardinality $\leq
\aleph _{1}$, then $\limfunc{Ext}(A,B^{*})=0.$
\end{theorem}

\smallskip\

We will make use of the identification of W-groups of cardinality
$\aleph_1$ as Shelah groups under the set-theoretic hypotheses. (See
\cite{Sh80} or \cite[XII.2.5, XIII.3.6]{EM}.) We fix a short exact
sequence
\[
0\rightarrow B^{*}\stackrel{\iota }{\hookrightarrow }N\stackrel{\pi }{%
\rightarrow }A\rightarrow 0
\]
and proceed to prove that it splits. We may assume that $\iota $ is
the inclusion map. We also fix a set function $\gamma :A\rightarrow
N$ such that for all $a\in A$, $\pi (\gamma (a))=a$. We will show
that the short exact sequence splits by proving the existence of a
function $h:A\rightarrow B^{*}$ such that the function
\[
\gamma -h:A\rightarrow N:a\mapsto \gamma (a)-h(a)
\]
is a homomorphism. The function $h$ will be obtained via a directed subset $%
\mathcal{G}$ of a c.c.c. poset $P$; the directed subset (which will
be required to intersect certain dense subsets of $P$) will exist as
a consequence of MA.

We define $P$ to consist of all triples $p=(A_{p},B_{p},h_{p})$
where $A_{p}$ (resp. $B_{p}$) is a pure and finitely generated
summand of $A$ (resp. $B$) and $h_{p}$ is a function from $A_{p}$ to
$B_{p}^{*}$. Moreover, we require that the function which takes
$a\in A_{p}$ to $\gamma (a)-h_{p}(a)$ is a homomorphism from $A_{p}$
into $N/\{f\in B^{*}:f\upharpoonright B_{p}\equiv 0\}$. (Strictly
speaking, this is an abuse of notation: by $\gamma (a)-h_{p}(a)$ we
mean the coset of $\gamma (a)-\eta _{a}$ where $\eta _{a}$
is any element of $B^{*}$ such that $\eta _{a}\upharpoonright B_{p}=h_{p}(a)$%
.)

The partial ordering on $P$ is defined as follows:
$p=(A_{p},B_{p},h_{p})\leq p^{\prime }=(A_{p}^{\prime
},B_{p}^{\prime },h_{p}^{\prime })$ if and only
if $A_{p}\subseteq A_{p}^{\prime }$, $B_{p}\subseteq B_{p}^{\prime }$, and $%
h_{p}^{\prime }(a)\upharpoonright B_{p}=h_{p}(a)$ for all $a\in
A_{p}$. The dense subsets that we use are:

\[
D_{a}^{1}=\{p\in P:a\in A_{p}\}
\]
for all $a\in A$ and
\[
D_{b}^{2}=\{p\in P:b\in B_{p}\}
\]
for all $b\in B$. Assuming that these sets are dense and that $P$ is
c.c.c., the axiom MA + $\lnot $CH yields a directed subset
$\mathcal{G}$ which has
non-empty intersection with each of these dense subsets. We can then define $h$ by: $%
h(a)(b)=h_{p}(a)(b)$ for some (all) $p\in \mathcal{G}$ such that
$a\in A_{p}$ and $b\in B_{p}$. It is easy to check that $h$ is
well-defined and has the desired properties.

For use in proving both density and the c.c.c. property, we state
the following claim, whose proof we defer to the end.

\begin{quotation}
(1) \textit{Given a basis }$\{x_{i}:i=1,...,n\}$ \textit{of a
finitely generated pure subgroup }$A_{0}$\textit{\ of }$A$\textit{, a basis }%
$\{y_{j}:j=1,...,m\}$\textit{\ of a finitely generated pure subgroup }$B_{0}$%
\textit{\ of }$B$\textit{, and an indexed set }$\{e_{ij}:i=1,...,n$, $%
j=1,...,m\}$\textit{\ of elements of }$\mathbb{Z}$\textit{, there is one and only one }$%
p\in P$\textit{\ such that }$A_{p}=A_{0}$\textit{,
}$B_{p}=B_{0}$\textit{\
and }$h_{p}(x_{i})(y_{j})=e_{ij}$\textit{\ for all }$i=1,...,n$\textit{, }$%
j=1,...,m$\textit{.}
\end{quotation}

Assuming this claim, we proceed to prove the density of $D_{a}^{1}$. Given $%
p $, there is a finitely generated pure subgroup $A^{\prime }$ of
$A$ which contains $A_{p}$ and $a$. Now $A_p$ is a summand of
$A^{\prime}$ so we can choose a basis $\{x_{i}:i=1,...,n\}$ of
$A^{\prime
}$ which includes a basis $\{x_{i}:i=1,...,k\}$ of $A_{p}$; choose a basis $%
\{y_{j}:j=1,...,m\}$ of $B_{p}$. Then by the claim there is an element $%
p^{\prime }$ of $P$ such that $A_{p^{\prime }}=A^{\prime }$,
$B_{p^{\prime
}}=B_{p}$, and $h_{p^{\prime }}(x_{i})(y_{j})=h_{p}(x_{i})(y_{j})$ for all $%
i=1,...,k$. Clearly $p^{\prime }$ extends $p$ and belongs to
$D_{a}^{1}$. The proof of the density of $D_{b}^{2}$ is similar.

Next we prove that $P$ is c.c.c. We will make use of the following
fact, which is proved in \cite{E}, Lemma 7.5. (It is proved there
for strongly $\aleph_1$-free groups, there called groups with
Chase's condition, but the proof may be adapted for Shelah groups;
cf. \cite[Theorem 7.1]{E2}.)

\begin{quotation}
(2) \textit{If }$G$\textit{\ is a Shelah group of cardinality }$\aleph _{1}$%
\textit{\ and }$\{S_{\alpha }:\alpha \in \omega _{1}\}$\textit{\ is
a family of finitely generated pure subgroups of }$G$\textit{, then
there is an uncountable subset }$I$\textit{\ of }$\omega
_{1}$\textit{\ and a pure free subgroup }$G^{\prime }$ of $G$
\textit{such that }$S_{\alpha }\subseteq G^{\prime }$\textit{\ for
all }$\alpha \in I$\textit{.}
\end{quotation}

Suppose that $\{p_{\nu }=(A_{\nu ,}B_{\nu },h_{\nu }):\nu \in \omega
_{1}\}$
is a subset of $P$. We must prove that there are indices $\mu \neq \nu $ such that $%
p_{\mu }$ and $p_{\nu }$ are compatible. By claim (2), passing to a
subset, we can assume that there is a pure free subgroup $A^{\prime
}$ of $A$ and a pure free subgroup $B^{\prime }$ of $B$ such that
$A_{\nu }\subseteq A^{\prime }$ and $B_{\nu }\subseteq B^{\prime }$
for all $\nu \in \omega _{1} $. (Now we follow the argument in
\cite{E} for property (7.1.3).) Choose a basis $X$ of $A^{\prime }$
and a basis $Y$ of $B^{\prime }$. By density we can assume that each
$A_{\nu }$ is generated by a finite subset $X_{\nu} $ of $X$ and
each $B_{\nu }$ is generated by a finite subset, $Y_{\nu }$, of $Y$.
Moreover we can assume that there is a (finite) subset $T$ of $X$
(resp. $W$ of $Y$) which is contained in each $X_{\nu }$ (resp. each
$Y_{\nu }$) and is maximal with respect to the property that it is
contained in uncountably many $X_{\nu }$ (resp. $Y_{\nu }$). Passing
to a subset, we can assume that $h_{\nu }(x)(y)$ has a value
independent of $\nu $ for each $x\in T$ and $y \in W$. By a counting
argument we can find $\nu
>0$ such that $X_{\nu }\cap X_{0}=T$ and $Y_{\nu }\cap Y_{0}=W$. We
define a
member $q\in P$ which extends $p_{0}$ and $p_{\nu }$, as follows. Let $%
A_{q}=\left\langle X_{0}\cup X_{\nu }\right\rangle $ and
$B_{q}=\left\langle Y_{0}\cup Y_{\nu }\right\rangle $. Clearly these
are pure subgroups of $A$
(resp. $B$). Then by claim (1) there is a homomorphism $h_{q}:A_{q}%
\rightarrow B_{q}^{*}$ such that $q\geq p_{0}$ and $q\geq p_{\nu }$.
In
particular, for $x\in T$%
\[
h_{q}(x)(y)=\left\{
\begin{tabular}{ll}
the common value & if $y\in W$ \\
$h_{0}(x)(y)$ & if $y\in Y_{0}-W$ \\
$h_{\nu }(x)(y)$ & if $y\in Y_{\nu }-W$%
\end{tabular}
\right.
\]
and for $x\in X_{0}-T$%
\[
h_{q}(x)(y)=\left\{
\begin{tabular}{ll}
$h_{0}(x)(y)$ & if $y\in Y_{0}$ \\
arbitrary & if $y\in Y_{\nu }-W$%
\end{tabular}
\right.
\]
and similarly for $x\in X_{\nu }-T$.

Thus we have shown the existence of $p_{\mu }$ and $p_{\nu }$ which
are compatible, and it remains to prove claim (1). We will prove
uniqueness first. Suppose that, with the notation of (1), there are
$p_{1}$ and $p_{2}$ in $P$ such that for $\ell =1,2$, $h_{p_{\ell
}}(x_{i})(y_{j})=e_{ij}$
 for all $i=1,...,n$\, $j=1,...,m$. It
suffices to prove that for all $a\in A_{0}$,
\[
h_{p_{1}}(a)(y_{j})=h_{p_{2}}(a)(y_{j})
\]
for all $j=1,...,m$. Write $a$ as a linear combination of the basis: $%
a=\sum_{i=1}^{n}d_{i}x_{i}$. By the definition of $P$%
\[
\gamma (a)-\sum_{i=1}^{n}d_{i}\gamma (x_{i})=h_{p_{\ell
}}(a)-\sum_{i=1}^{n}d_{i}h_{p_{\ell }}(x_{i}).
\]
for $\ell =1,2$. Note that the right-hand side is, by definition,  a
function on $B_0$.
Applying both sides to $%
y_{j}\in B_0$, we obtain that
\[
h_{p_{\ell }}(a)(y_{j})=\sum_{i=1}^{n}d_{i}e_{ij}+(\gamma
(a)-\sum_{i=1}^{n}d_{i}\gamma (x_{i}))(y_{j})
\]
for $\ell =1,2$. Since the right-hand side is independent of $\ell
$, we can conclude the desired identity. For existence, we use the
last displayed equation to define $h_{p}$ and easily check that it
gives an element of $P$.  \qed

\medskip

It is possible to force to obtain

\noindent \textbf{Model (3).} A model of MA + $2^{\aleph_{0}} =
\aleph_2$ plus

\begin{quote}
 ($\ddagger $) $\diamondsuit _{\kappa }(E)$ holds for every
regular cardinal $\kappa >\aleph _{1}$ and every stationary subset $E$ of $%
\kappa $.
\end{quote}

\medskip

Just as for Model (2), in this model the tensor product of any two
Whitehead groups of arbitrary cardinality is again a W-group.


\begin{thebibliography}{99}
\bibitem{BFS} T. Becker, L.  Fuchs,  and S. Shelah,
 {\it Whitehead modules over domains}, Forum Math. {\bf 1} (1989), 53--68.

\bibitem{CE} H. Cartan,  and S. Eilenberg,    {\bf Homological Algebra},
Princeton University Press (1956).

\bibitem{E} P. C. Eklof, {\it Whitehead's problem is undecidable},
Amer. Math. Monthly {\bf 83} (1976), 775--788.

\bibitem{E2} P. C. Eklof, {\bf Set Theoretic Methods in Homological
Algebra and Abelian Groups}, Les Presses de L'Universit\'e de
Montr\'eal (1980).

\bibitem{EM}  P. C. Eklof and A. H. Mekler, \textbf{Almost Free Modules}, Rev. Ed,
North-Holland (2002).

\bibitem{ES}  P.C. Eklof and S. Shelah, \textit{The structure of} Ext\textit{$(A, \mathbb{Z})$ and
GCH: possible co-Moore spaces}, Math. Zeit. \textbf{239} (2002),
143--157.

\bibitem{F}  L. Fuchs, \textbf{Infinite Abelian Groups}, vols I and II,
Academic Press (1970, 1973).




\bibitem{Hu}  M. Huber, \textit{On reflexive modules and abelian groups}, J.
Algebra \textbf{82} (1983), 469--487.





\bibitem{Sh74}  S. Shelah, \textit{Infinite abelian groups, Whitehead
problem and some constructions}, Israel J. Math \textbf{18} (1974),
243--256.





\bibitem{Sh80}  S.{\ Shelah, }\textit{Whitehead groups may not be free even
assuming CH, II}, Israel J. Math. \textbf{35 (}1980\textbf{),} 257--285.


\end{thebibliography}
\end{document}